\documentclass[12pt]{amsart}
\usepackage{graphicx, overpic, epsf,amssymb,amscd,amsfonts,times}

\newtheorem{thm}{Theorem}[section]

\newtheorem{prop}[thm]{Proposition}

\newtheorem{lemma}[thm]{Lemma}

\newcommand{\C}{\mathbb{C}}
\newcommand{\R}{\mathbb{R}}
\newcommand{\Z}{\mathbb{Z}}

\newcommand{\bdry}{\partial}

\newcommand{\s}{\vskip.1in}
\newcommand{\n}{\noindent}

\newcommand{\be}{\begin{enumerate}}
\newcommand{\ee}{\end{enumerate}}
\newcommand{\veer}{Veer(S,\bdry S)}

\newcommand{\spinc}{$\mbox{Spin}^c$ \hskip-.02in}
\newcommand{\T}{\mathbb{T}}

\numberwithin{equation}{subsection}

\begin{document}

\title{On the contact class in Heegaard Floer homology}

\author{Ko Honda}
\address{University of Southern California, Los Angeles, CA 90089}
\email{khonda@usc.edu} \urladdr{http://rcf.usc.edu/\char126 khonda}

\author{William H. Kazez}
\address{University of Georgia, Athens, GA 30602} \email{will@math.uga.edu}
\urladdr{http://www.math.uga.edu/\char126 will}

\author{Gordana Mati\'c}
\address{University of Georgia, Athens, GA 30602} \email{gordana@math.uga.edu}
\urladdr{http://www.math.uga.edu/\char126 gordana}

\date{This version: March 6, 2007. (The pictures are in color.)}



\thanks{KH supported by an Alfred P.\ Sloan Fellowship and an NSF
CAREER Award (DMS-0237386); GM supported by NSF grant DMS-0410066;
WHK supported by NSF grant DMS-0406158.}

\begin{abstract}
We present an alternate description of the Ozsv\'ath-Szab\'o contact
class in Heegaard Floer homology.  Using our contact class, we prove
that if a contact structure $(M,\xi)$ has an adapted open book
decomposition whose page $S$ is a once-punctured torus, then the
monodromy is right-veering if and only if the contact structure is
tight.
\end{abstract}

\maketitle

\section{Introduction}

In the paper \cite{OS5}, Ozsv\'ath and Szab\'o defined an invariant
of a contact 3-manifold $(M,\xi)$ which lives in the Heegaard Floer
homology $\widehat{HF}(-M)$ of the manifold $M$ with reversed
orientation.  It is defined via the work of Giroux~\cite{Gi2}, who
showed that there is a 1-1 correspondence between isomorphism
classes of open book decompositions modulo positive stabilization
and isomorphism classes of contact structures on closed 3-manifolds.
Ozsv\'ath and Szab\'o associated an element in Heegaard Floer
homology to an open book decomposition  and showed that its homology
class is independent of the choice of the open book compatible with
the given contact structure. They also showed that this invariant
$c(\xi)$ is zero if the contact structure is overtwisted, and that
it is nonzero if the contact structure is symplectically fillable.
The {\em contact class} $c(\xi)$ has proven to be extremely powerful
at (i) proving the tightness of various contact structures and (ii)
distinguishing tight contact structures, especially in the hands of
Lisca-Stipsicz~\cite{LS1,LS2} and Ghiggini~\cite{Gh}.

The goal of this paper is to introduce an alternate, more hands-on,
description of the contact class in Heegaard Floer homology and to
use it in the context of our program of relating right-veering
diffeomorphisms to tight contact structures.

In \cite{HKM2} we introduced the study of right-veering
diffeomorphisms of a compact oriented surface with nonempty boundary
(sometimes called a ``bordered surface''), and proved that if
$(S,h)$ is an open book decomposition compatible with a tight
contact structure, then $h$ is right-veering.  In \cite{HKM3} we
continued the study of the monoid $\veer$ of right-veering
diffeomorphisms and investigated its relationship with symplectic
fillability in the pseudo-Anosov case. We proved the following:

\begin{thm} \label{cgeq1}
Let $S$ be a bordered surface with connected boundary and $h$ be
pseudo-Anosov with fractional Dehn twist coefficient $c$. If $c\geq
1$, then the contact structure $\xi_{(S,h)}$ supported by $(S,h)$ is
isotopic to a perturbation of a taut foliation. Hence $(S,h)$ is
(weakly) symplectically fillable and universally tight if $c\geq 1$.
\end{thm}

Hence, when a contact structure is supported by an open book with
``sufficiently'' right-veering monodromy, it is symplectically
fillable and therefore tight as a consequence of a theorem of
Eliashberg and Gromov~\cite{El}.  Unfortunately, a right-veering
diffeomorphism with a small amount of rotation does not always
correspond to a tight contact structure. In fact, any open book can
be stabilized to a right-veering one (see Goodman~\cite{Go}, as well
as \cite{HKM2}).  However, we might optimistically conjecture that a
minimal (i.e., not destabilizable) right-veering open book defines a
tight contact structure.  If we specialize to the case of a
once-punctured torus, then we can use our description of the contact
class to prove this conjecture.

\begin{thm} \label{torus}
Let $(M,\xi)$ be a contact 3-manifold which is supported by an open
book decomposition $(S,h)$, where $S$ is a once-punctured torus.
Then $\xi$ is tight if and only if $h$ is right-veering.
\end{thm}

Very recently John Baldwin~\cite{Ba} also obtained results similar
to Theorem~\ref{torus}.

The paper is organized as follows. In Section~\ref{review}, we
review the standard definition of $c(\xi)$.  Then, in
Section~\ref{alternate}, we describe the class $EH(\xi)\in
\widehat{HF}(-M)$, which arose in discussions between John Etnyre
and the first author.  We also prove that the class $EH(\xi)$ equals
the Ozsv\'ath-Szab\'o contact class $c(\xi)$, and hence $EH(\xi)$ is
a contact invariant. In Section~\ref{EH}, the class $EH(\xi)$ is
applied to contact structures with compatible genus one open book
decompositions to prove Theorem~\ref{torus}.

\section{Open books and Ozsv\'ath-Szab\'o contact invariants}
\label{review}

In \cite{OS1,OS2}, Ozsv\'ath and Szab\'o defined invariants of
closed oriented 3-manifolds $M$ which they called {\em Heegaard
Floer homology}.  Among the several versions of Heegaard Floer
homology defined by Ozsv\'ath and Szab\'o, we concentrate on the
simplest one, namely $\widehat {HF}(M)$. It is defined as the
homology associated to a chain complex determined by a Heegaard
decomposition of $M$. Consider a Heegaard decomposition
$(\Sigma,\alpha=\{\alpha_1,\dots,\alpha_g\},
\beta=\{\beta_1,\dots,\beta_g\})$ of genus $g$. Here $\Sigma$ is the
Heegaard surface, i.e., a closed oriented surface of genus $g$ which
splits $M$ into two handlebodies $H_1$ and $H_2$, $\Sigma=\bdry
H_1=-\bdry H_2$, $\alpha_i$ are the boundaries of the compressing
disks of $H_1$, and $\beta_i$ are the boundaries of the compressing
disks of $H_2$. Then consider two tori $\T_{\alpha}=\alpha_1 \times
\dots \times \alpha_g$ and $\T_{\beta}=\beta_1 \times \dots \times
\beta_g$ in $Sym^g(\Sigma)$. Also pick a basepoint $z\in \Sigma$.
The complex $\widehat{CF}(M)$ is defined to be the free $\Z$-module
generated by the points $\mathbf{x} = (x_1, \dots , x_g)$ of
$\T_{\alpha}\cap\T_{\beta}$. The boundary map is defined by counting
points in certain 0-dimensional moduli spaces of holomorphic maps of
the unit disk into $Sym^g(\Sigma)$. It is, very roughly, defined as
follows. Denote by $\mathcal{M}_{\mathbf{x} ,\mathbf{y}}$ the
0-dimensional (after quotienting by the natural $\R$-action) moduli
space of holomorphic maps $u$ from the unit disk $D^2\subset \C$ to
$Sym^g(\Sigma)$ that (i) send $1\mapsto \mathbf{x}$, $-1\mapsto
\mathbf{y}$, $S^1 \cap \{\mbox{Im } z \geq 0 \}$ to $\T_{\alpha}$
and $S^1 \cap \{\mbox{Im } z \leq 0\}$ to $\T_{\beta}$, and (ii)
avoid $\{z\}\times Sym^{g-1}(\Sigma)\subset Sym^{g}(\Sigma)$. Then
define
$$\bdry \mathbf{x} =  \sum_{\mu(\mathbf{x},\mathbf{y})=1} ~~
\#(\mathcal{M}_{\mathbf{x}, \mathbf{y} })~~  \mathbf{y},$$ where
$\mu(\mathbf{x},\mathbf{y})$ is the relative Maslov index of the
pair and $\#(\mathcal{M}_{\mathbf{x}, \mathbf{y} })$ is a signed
count of points in $\mathcal{M}_{\mathbf{x}, \mathbf{y} }$. The
homology of this complex $\widehat{HF}(M)$ is shown to be
independent of the various choices made in the definition. In
particular, it is independent of the choice of a ``weakly
admissible'' Heegaard decomposition.

Each intersection point $\mathbf{x}$  in $\T_\alpha\cap \T_\beta$
defines a \spinc structure $\mathbf{s}_\mathbf{x}$ on $M$. If there
is a topological disk from $\mathbf{x}$ to $\mathbf{y}$ which
satisfies (i) and (ii) in the previous paragraph, then the two
\spinc structures agree. Hence, the complex (as well as the homology
of the complex) splits according to \spinc structures. The Heegaard
Floer homology decomposes as a direct sum
$$\widehat {HF}(M) =  \oplus_{\mathbf{s}} ~~  \widehat {HF}(M, \bold s).$$

\begin{figure}[ht]
\s
\begin{overpic}[width=15cm]{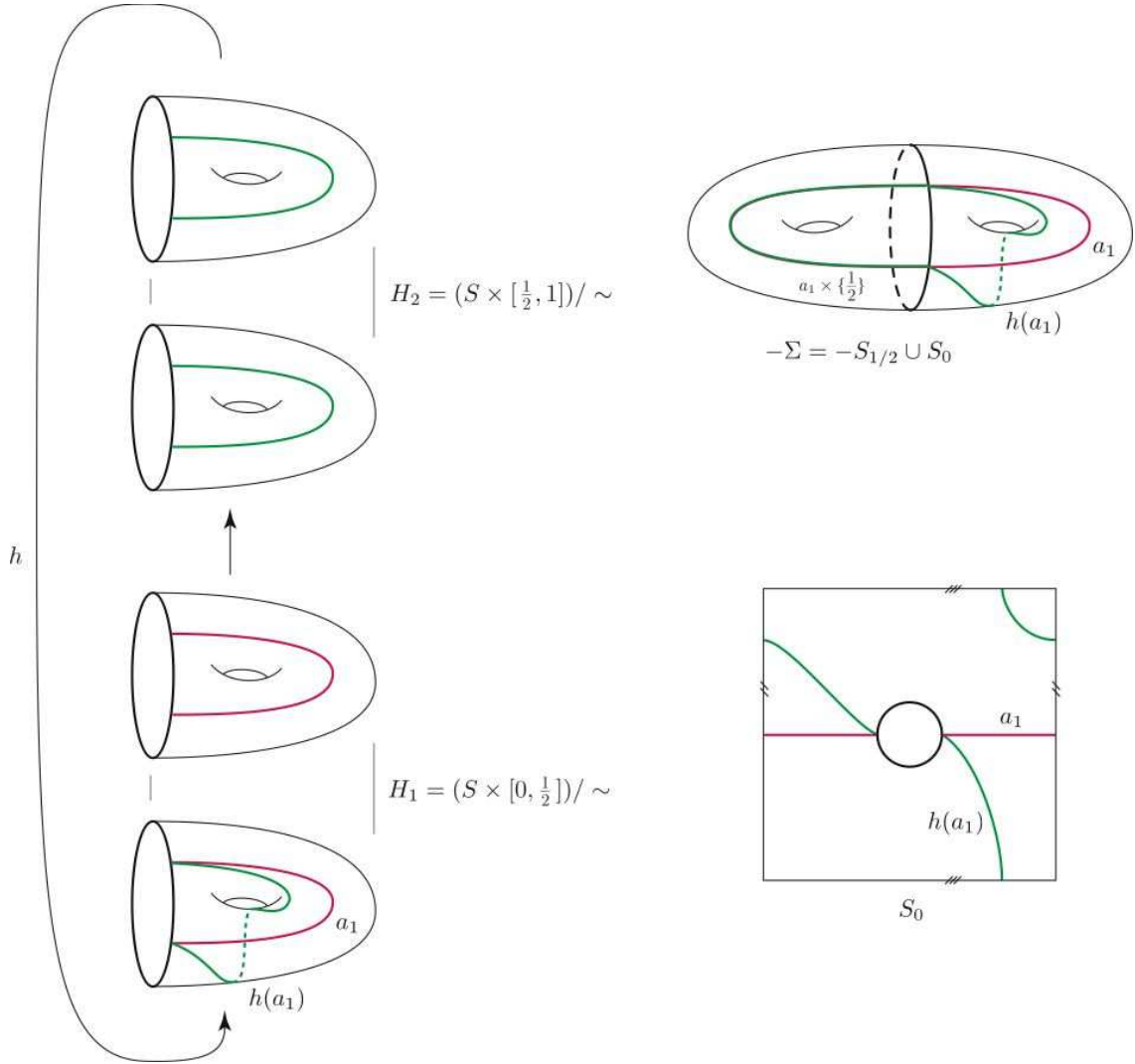}
\end{overpic}
\caption{The left-hand portion of the figure shows the decomposition
into the two handlebodies $H_1$ and $H_2$ and a compressing disk on
each corresponding to $a_1$. The upper right portion shows $-\Sigma
= -S_{1/2} \cup  S_0$ and the boundaries of two compressing disks.
We draw just the lower right portion to indicate the Heegaard
decomposition and the effect of the monodromy on arcs.}
\label{openbook}
\end{figure}

Given a contact structure $\xi$ on $M$, we denote the associated
\spinc structure by $\mathbf{s}_\xi$.  Let $(S,h,K)$ be an open book
decomposition of a manifold that is compatible with the contact
structure $\xi$.  Then Ozsv\'ath and Szab\'o define in \cite{OS5} an
element $c(\xi)\in \widehat {HF}(-M, \mathbf{s}_\xi)/(\pm 1)$ by
using a Heegaard splitting associated to the open book decomposition
as follows.  (At the time of the writing of the paper, the $\pm 1$
ambiguity still exists.  It is possible, however, that a careful
study of orientations would remove this ambiguity.  The $\pm 1$
issue does not arise in Seiberg-Witten Floer homology.) To avoid
writing $\pm 1$ everywhere, we either work with
$\Z/2\Z$-coefficients or tacitly assume that $c(\xi)$ is
well-defined up to a sign when $\Z$-coefficients are used. Consider
the open book decomposition $(S,h,K)$, where $S$ is a surface of
genus $g$ (here genus means the genus of the surface capped off with
disks) with one boundary component $\partial S$, $h$ is a
diffeomorphism of $S$ which is the identity on $\bdry S$, and the
pair $(M,K)$ is homeomorphic to $((S \times[0,1]) /\sim, (\partial S
\times[0,1]) /\sim)$. The equivalence relation $\sim$ is generated
by $(x,1) \sim (h(x),0)$ for $x\in S$ and $(y,t) \sim (y,t')$ for $y
\in \partial S$, $t,t'\in[0,1]$.  From the above description of $M$
we immediately see an associated Heegaard splitting of $M$ by
setting $H_1 =  (S \times[0,\frac{1}{2}]) /\sim$ and $H_2 = (S
\times[{1\over 2},1]) /\sim$. This gives a Heegaard decomposition of
genus $2g$ with the splitting surface $\Sigma = S_{1/2} \cup  -S_0$.
A set of $2g$ properly embedded disjoint arcs $a_1, \dots, a_{2g}$
which cut $S$ into a disk defines a set of compressing disks $a_i
\times [0,{1 \over 2}]$, $i=1,\dots,2g$, in $H_1$ and a set of
compressing disks $a_i\times [{1\over 2},1]$, $i=1,\dots,2g$, in
$H_2$. We then set $\alpha_i=\bdry (a_i\times[0,{1\over 2}])$ and
$\beta_i=\bdry (a_i\times [{1\over 2},1])$, for $i=1,\dots,2g$. See
Figure~\ref{openbook}.

\s This is, however, not the Heegaard splitting that Ozsv\'ath and
Szab\'o consider when defining $c(\xi)$.  Instead they use a
Heegaard surface that can be viewed simultaneously as a Heegaard
surface for $M$ and for $M_0(K)$, the zero surgery along the binding
$K$. The contact element in $\widehat{HF}(-M)$ can be seen on this
Heegaard surface as the image of a class in $\widehat{HF}(-M_0(K))$
(or, equivalently, as the image of a class in
$\widehat{HFK}(-M,K,F,-g)$). To construct such a splitting, take a
disk $D\subset int(S)$ which is contained in a small neighborhood of
$\bdry S$, dig $D \times [0,{1\over 2}]$ out of $H_1$, and then
attach it to $H_2$.  This produces two new handlebodies $H'_1$ and
$H'_2$. On $H'_2$ we keep the same set of $\beta$-curves
$\beta_1,\dots,\beta_{2g}$ as $H_2$ and add $\beta_0=\bdry D\times
\{{1\over 4}\}$. Next, let $d$ be a short arc connecting between the
two boundary components of $S-D$, and let $\{b_1,\dots,b_{2g}\}$ be
a set of arcs with endpoints on $\bdry D$ which are ``dual'' to
$\{a_1,\dots,a_{2g}\}$. (By this we mean $a_{2i+1}\cap
b_j=\emptyset$ if $j\not=2i$ and $a_{2i+1}\cap b_{2i}=\{x_{2i+1}\}$;
also $a_{2i}\cap b_j=\emptyset$ if $j\not=2i+1$ and $a_{2i}\cap
b_{2i+1}=\{x_{2i}\}$.) Then on $H'_1$, we let $\alpha_0=\bdry
(d\times[0,{1\over 2}])$ and $\alpha_i=\bdry (b_i\times [0,{1\over
2}])$. Also let $\alpha_0\cap \beta_0= \{x_0\}$.

These above choices determine a special point $\mathbf{x}=(x_0,x_1,
\dots,x_{2g})$ in $\T_\alpha\cap \T_\beta \subset
Sym^{2g+1}(\Sigma)$. (Here, $x_i$ means $(x_i, {1\over 2})$, for
$i>0$.) This point (after modifying the Heegaard diagram by winding
in a region that does not affect $\bf x$ to adjust for
admissibility) defines the special cycle in Heegaard Floer homology.
The homology class of $\bf x$ is defined as the contact class
$c(\xi)$ by Ozsv\'ath-Szab\'o.  They show that
$\widehat{HFK}(-M,K,F,-g)$, the knot Floer homology for $(-M,K)$ at
the lowest possible filtration level $-g$, is isomorphic to $\mathbb
Z$ and is generated by $\bf x$.  Then $c(\xi)$ is defined to be the
image of this generator in $\widehat{HF}(-M)$. For details,
including the figures describing this decomposition and the
corresponding generator of $c(\xi)$, see \cite{OS5}.

\section{An alternate description of the contact element}
\label{alternate}

\subsection{Definition and main theorem}
\begin{figure}[ht]
\begin{overpic}[width=5cm]{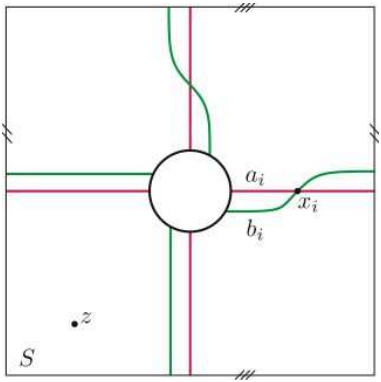}
\end{overpic}
\caption{The arcs $a_i$ and $b_i$ for a once-punctured torus $S$.}
\label{pushoff}
\end{figure}

Let $S$ be a bordered surface whose boundary is not necessarily
connected. Let $\{a_1,\dots,a_r\}$ be a collection of disjoint,
properly embedded arcs of $S$ so that $S-\bigcup_{i=1}^r a_i$ is a
single polygon. We will call such a collection a {\em basis for
$S$}. Observe that every arc $a_i$ of a basis is a nonseparating arc
of $S$. Next let $b_i$ be an arc which is isotopic to $a_i$ by a
small isotopy so that the following hold: \be
\item The endpoints of $a_i$ are isotoped along $\bdry S$, in the
direction given by the boundary orientation of $S$.
\item $a_i$ and $b_i$ intersect transversely in one point in the
interior of $S$.
\item If we orient $a_i$, and $b_i$ is given the induced orientation
from the isotopy, then the sign of the intersection $a_i\cap b_i$ is
$+1$. \ee See Figure~\ref{pushoff}.

Let $M=M_{(S,h)}$ be the 3-manifold with open book decomposition
$(S,h)$.  Recall the Heegaard decomposition for $M$ described in the
previous section, where $\Sigma=S_{1/2} \cup -S_0$.  We choose the
compressing disks to be $\alpha_i = \bdry (a_i \times[0,{1\over
2}])$ and $\beta_i = \bdry( b_i \times[{1\over 2},1])$.  We will
sometimes write $\alpha_i=(a_i,a_i)$ and $\beta_i=(b_i,h(b_i))$,
where the first entry is the arc on $S_{1/2}$ and the second entry
is the arc on $S_0$. Let $x_i$ be the intersection point $(a_i\cap
b_i)\times \{{1\over 2}\}$ lying in $S_{1/2} \subset \Sigma$, and
let $z$ be the basepoint which sits on $S_{1/2}$ and lies {\em
outside} the thin strips of isotopy between the $a_i$'s and the
$b_i$'s. Then $(\Sigma,\beta,\alpha,z)$ gives a weakly admissible
Heegaard diagram, namely every periodic domain has positive and
negative components. This is due to the fact that every periodic
domain which involves $\alpha_i$ crosses $x_i$, at which point the
sign of the connected component of $\Sigma-\bigcup_{i=1}^r
\alpha_i-\bigcup_{i=1}^r\beta_i$ changes.

Throughout this paper we use a product complex structure
$J=Sym^r(\mathfrak{j})$ on $Sym^r(\Sigma)$, where $\mathfrak{j}$ is
some complex structure on $\Sigma$, and perturb the $\alpha$- and
$\beta$-curves to attain transversality.  This is done using the
technique of \cite{Oh}, as sketched in Section~3.5 of \cite{OS1}. We
remark that moving the $\alpha$- and $\beta$-curves represents a
subclass of the Hamiltonian isotopies of $\T_\alpha$ and $\T_\beta$
(i.e., we have fewer perturbations), so Theorem~I of \cite{Oh} does
not carry over verbatim, but the proof technique carries over
without difficulty. Observe that if there is no holomorphic disk in
a given homotopy class, then the moduli space of such disks is
automatically Fredholm regular.

A $J$-holomorphic disk $u:D\rightarrow Sym^r(\Sigma)$ corresponds to
a holomorphic map $\hat u: \widehat D\rightarrow \Sigma$, where
$\widehat D$ is a branched cover of $D$.  In the definition of the
boundary map in the $\widehat{HF}$ theory, we only count holomorphic
disks $u:D\rightarrow Sym^r(\Sigma)$ that miss $\{z\}\times
Sym^{r-1}(\Sigma)$.  Hence it follows that we only count $\hat u$
for which the image of $\hat u$ does not intersect $z\in \Sigma$.
The intersection of any such $\hat u$ with $S_{1/2}$ is thus
constrained to lie in the thin strips of isotopy of the $a_i$ to
$b_i$.

We claim that $\mathbf{x}=(x_1,\dots,x_r)\in
\widehat{CF}(\Sigma,\beta,\alpha,z)$ is a cycle, thanks to the
fortuitous placement of the basepoint $z$. (We write
$\widehat{CF}(\Sigma,\beta,\alpha,z)$ instead of
$\widehat{CF}(\Sigma,\alpha,\beta,z)$ to indicate homology on $-M$.)
Suppose $\hat u$ contributes to $\bdry \mathbf{x}$; in particular it
is nonconstant.  Let $\delta_i$ be a short oriented arc of $\bdry
\widehat D$ which passes through a corner $p_i\in\widehat D$ for
which $\hat{u}(p_i)= x_i$.  Then $\hat u(\delta_i)$ first travels
along $\alpha_i$ and switches to $\beta_i$ at $x_i$. More
explicitly, there is some $t_0$, $\delta_i(t_0)\in \alpha_i$, such
that ${d\over dt}(\hat{u}\circ \delta_i)(t_0)\not=0$ and points
towards $x_i$. Since the interior of $\widehat D$ is to the left of
$\delta_i$, by the openness of the holomorphic map, $z$ would be
contained in the image of $\hat u$, a contradiction.

We define $EH(S,h,\{a_1,\dots,a_r\})$ to be the homology class of
the generator $\mathbf{x}$. The following is the main theorem of
this section:

\begin{thm}\label{EH=os}
$EH(S,h,\{a_1,\dots,a_r\})$ is an invariant of the contact structure
and equals $c(\xi_{(S,h)})$.
\end{thm}

In particular, $EH(S,h,\{a_1,\dots,a_r\})$ is independent of the
choice of basis, and it will often be denoted by $EH(S,h)$.

In Theorem~\ref{EH=os} we are not assuming that $\bdry S$ is connected.

\begin{figure}[ht]
\begin{overpic}[width=12cm]{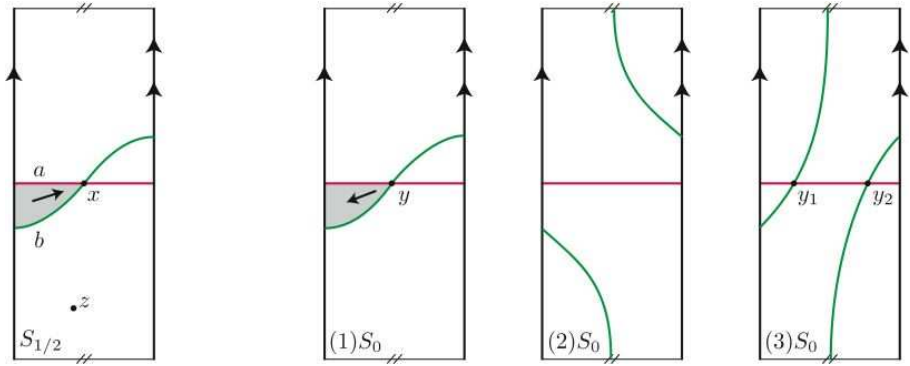}
\end{overpic}
\caption{Examples when $S$ is an annulus.} \label{examples}
\end{figure}

\s\n {\bf Examples:} To give some intuition for the class $EH(S,h)$, we
give three examples when $S$ is an annulus.  Refer to
Figure~\ref{examples}.  The leftmost diagram gives $a$ and $b$ on
$S_{1/2}$.  The subsequent diagrams give $S_{0}$ for (1), (2), and
(3) below (from left to right).

\be

\item If $h$ is the identity, then $(M,\xi)$ is the standard tight
contact structure on $S^1\times S^2$. Since there are two
holomorphic disks from $y$ to $x$, it follows that $EH(S,h)\not=0$.
One of the holomorphic disks from $y$ to $x$ has been shaded in
Figure~\ref{examples}.

\item If $h$ is a positive Dehn twist about the core curve, then
$(M,\xi)$ is the standard tight contact structure on $S^3$. Since
$x$ is the unique intersection point on $\Sigma=T^2$,
$EH(S,h)\not=0$.

\item If $h$ is a negative Dehn twist about the core curve, then
$(M,\xi)$ is an overtwisted contact structure on $S^3$. We have
$\bdry y_1=\bdry y_2=x$; hence $EH(S,h)=0$.

\ee

The following lemma echoes our result in \cite{HKM2}, which states
that $\xi_{(S,h)}$ is overtwisted if $h$ is not right-veering.

\begin{lemma}\label{rv}
If $h$ is not right-veering, then $EH(S,h)=0$.
\end{lemma}

\begin{proof}
If $h$ is not right-veering, then there exists an arc $a_1$ on $S$
so that $h(a_1)$ is to the left of $a_1$.   If $a_1$ is
nonseparating, then it can be completed to a basis
$\{a_1,\dots,a_r\}$.  There exists an intersection point $y_1\in
\alpha_1\cap \beta_1$ and a unique (up to translation) holomorphic
disk $D\subset \Sigma$ from $y_1$ to $x_1$, where $1\mapsto y_1$,
$-1\mapsto x_1$, $\bdry D\cap \{y\geq 0\}$ maps to $\beta_i$ and
$\bdry D\cap \{y\leq 0\}$ maps to $\alpha_i$. Since $z$ forces any
holomorphic disk $\hat{u}:\widehat{D}\rightarrow \Sigma$ which
contributes to $\bdry (y_1,x_2,\dots, x_n)$ to be constant near
$x_i$, $i=2,\dots,r$, all the $\alpha_i$ and $\beta_i$,
$i=2,\dots,r$, are ``used up'', and the only holomorphic disk that
remains is the unique one from $y_1$ to $x_1$. Hence $\bdry
(y_1,x_2,\dots,x_n)=(x_1,x_2,\dots,x_n)$.

If the arc $a_1$ is separating, then let us call its initial point
$p$. The arcs $h(a_1)$ and $a_1$ must intersect at some point $q$ in
the interior of $a_1$; otherwise $h(a_1)$ will cut off a strictly
smaller subsurface of $S$ inside a subsurface of $S$ cut off by
$a_1$. Let $c$ be the subarc of $a_1$ from $p$ to $q$ and $c'$ be
the subarc of $h(a_1)$ from $p$ to $q$. Then either $c(c')^{-1}$ is
separating or it is not.  If $c(c')^{-1}$ separates a region $S'$ to
the left of $a_1$, then there is a nonseparating arc $b\subset S'$
which begins and ends at $p$. On the other hand, if $c(c')^{-1}$ is
nonseparating, then we let $b=c(c')^{-1}$.  In either case, since
$b$ is strictly to the left of $a_1$ and strictly to the right of
$h(a_1)$, it follows that $h(b)$ is strictly to the left of $b$.
\end{proof}

In view of Theorem~\ref{EH=os} and the fact that every overtwisted
contact structure admits an open book that is not right-veering,
Lemma~\ref{rv} immediately implies that $c(\xi)=0$ for an
overtwisted contact structure.

\begin{proof}[Proof of Theorem~\ref{EH=os}.]
Let us denote a positive Dehn twist about a closed curve $\gamma$ by
$\phi_\gamma$. Assume $\bdry S$ is connected. We first prove the
theorem for a special case, namely when $h=\phi_{\bdry S}^n$ with
$n>0$, in Section~\ref{primordial}. Next, in Section~\ref{isotopy}
we prove that $EH(S,h,\{a_1,\dots,a_r\})$ only depends on the
isotopy class of $h$ (relative to the boundary), and in
Section~\ref{change of basis} we show that $EH(S, h,\{a_1, \dots,
a_r\})$ is independent of the choice of basis by using handleslides.
Then in Section~\ref{Legendrian} we prove that $EH(S,h)$ is mapped
to $EH(S,\phi^{-1}_\gamma\circ h)$ under the natural map
$\widehat{HF}(-M_{(S,h)})\rightarrow
\widehat{HF}(-M_{(S,\phi^{-1}_\gamma\circ h)})$ which corresponds to
a Legendrian $(+1)$-surgery.  We then start with $\phi_{\bdry S}^n$
with $n\gg 0$ and apply a sequence of negative Dehn twists until we
reach the desired monodromy map $h$.  In Section~\ref{multiple} we
reduce the case of multiple boundary components to the case when
$\bdry S$ is connected.
\end{proof}

\subsection{Primordial Example} \label{primordial}

Let $S$ be a once-punctured torus and $h=\phi_{\bdry S}$, i.e., a
positive Dehn twist about $\bdry S$.  The same argument works if $S$
is a genus $g$ surface with one boundary component and
$h=\phi_{\bdry S}^n$, $n>0$.

\begin{figure}[ht]
\begin{overpic}[width=6in]{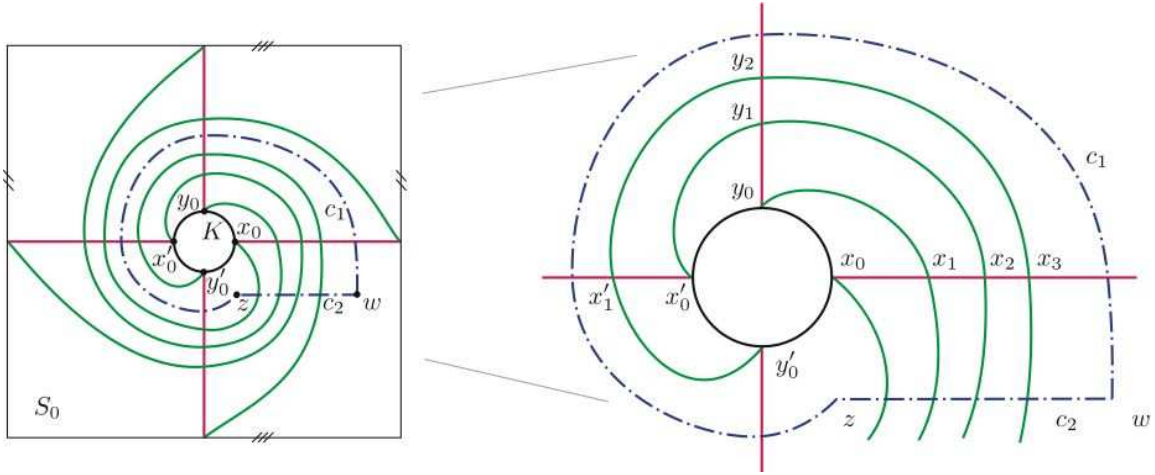}
\end{overpic}
\caption{$S_0$ when $h=\phi_{\bdry S}$, and a zoomed-in
ammonite-like region.} \label{twisted}
\end{figure}

The subarcs of $\alpha_i$ and $\beta_i$ that live in $S_0$ are given
in Figure~\ref{twisted}.  We change notation and the constituent
points of $\mathbf{x}$ representing $EH(S,h)$ will be denoted
$x_0=x_0'$ and $y_0=y_0'$ as in Figure~\ref{twisted}.  Although
$x_0=x_0'$ and $y_0=y_0'$, strictly speaking, live on $S_{1/2}$, we
view them as sitting on $\bdry S_0$. (Also, the points $x_0$ and
$x_0'$, as well as $y_0$ and $y_0'$, are drawn as distinct points on
$\bdry S_0$, but we hope this will not cause any confusion for the
reader.)

We then place the basepoint $w$ on $S_0$ as indicated in
Figure~\ref{twisted}. Observe that $z$ and $w$ together represent
the binding $K$.  The binding $K$ is isotopic to the dotted curve
$\gamma_0$ which consists of two subarcs $c_1$ and $c_2$ between $z$
and $w$, where $c_1$ intersects only $\alpha$-curves and $c_2$
intersects only $\beta$-curves.  Then $(\Sigma, \beta, \alpha, z,
w)$ is a doubly-pointed Heegaard diagram for the knot Floer homology
of $K$.

\begin{figure}[ht]
\begin{overpic}[width=6in]{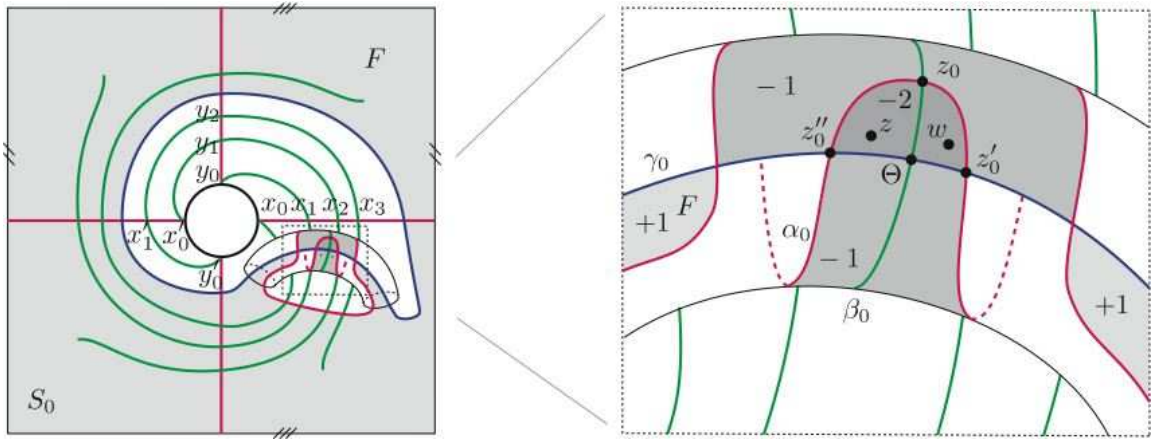}
\end{overpic}
\caption{Part of the stabilized Heegaard surface $-\Sigma'$. The
domain $F$ has been shaded.} \label{spiral}
\end{figure}

If we stabilize this Heegaard splitting by digging a handle in
$S\times[0,{1\over 2}]$ which is parallel to the arc $c_2$, then we
obtain a Heegaard surface $\Sigma'$ on which we can see both $-M$
and $-M_0(K)$. See Figure~\ref{spiral}. Here $-M$ is given by
$\{\beta_0\}\cup\beta$ and $\{\alpha_0\}\cup\alpha$, whereas
$-M_0(K)$ is given by $\gamma=\{\gamma_0\}\cup \beta$ and
$\{\alpha_0\}\cup\alpha$. (Here $\gamma_0$ is viewed as a curve that
passes through the handle once.) The stabilization sends
$\mathbf{x}= (x_0,y_0)$ to $\mathbf{x}'= (z_0, x_0,y_0)$, where
$z_0$ is the intersection of the two new compressing curves
$\alpha_0$ and $\beta_0$.

As a first step in exploiting the Ozsv\'ath-Szab\'o characterization
of $c(\xi)$, we show that the lowest filtration level is generated
by $\mathbf{x}'=(z_0,x_0,y_0)$ as well as the other intersection
points $\mathbf{y}=(z_0,x,y)$, where $x$ and $y$ live inside the
dotted lines of Figure~\ref{twisted}.  The filtration level is
computed by first letting $F\subset \Sigma'$ be the domain bounded
by $\alpha_0$ and $\gamma_0$ which does not intersect $S_{1/2}$ (and
hence lives mostly on $S_0$). We additionally assume that $F$ is
oriented so that the surface $\hat F$, obtained from $F$ by capping
off $\bdry F$, is an oriented fiber of the fibration of $M_0(K)$. In
order to find generators $\mathbf{y}$ which are at the lowest
filtration level, we minimize $\langle
c_1(\mathbf{s}_{\mathbf{y}'}),[\hat F]\rangle$. Here
$\mathbf{y}'=(z_0', x, y)$ and $z_0'$ is the point on
$\alpha_0\cap\gamma_0$ which is close to $z_0$ and obtained by
tensoring $z_0$ with the unique intersection point $\Theta\in
\beta_0\cap \gamma_0$ as in Figure~\ref{spiral}. (Keep in mind that
since we are dealing with $\widehat{HF}$ of $-M$ and $-M_0(K)$, the
Heegaard surface is $-\Sigma'$; otherwise our calculations will be
off by a negative sign.)

To this end, we recall the first Chern class formula (Section~7.1 of
\cite{OS2}; for some details, see Rasmussen~\cite{Ra}):
$$\langle c_1(\mathbf{s}_{\mathbf{y}'}),[A]\rangle=\chi(\mathcal{P})-2\overline{n}_{z}
(\mathcal{P}) +2\sum_{p\in \mathbf{y}'}
\overline{n}_{p}(\mathcal{P}).$$ Here $[A]\in H_2(M_0(K),\Z)$,
$\mathbf{s}_{\mathbf{y}'}$ is a \spinc structure corresponding to
$\mathbf{y}'$, $\mathcal{P}$ is the periodic domain for $[A]$ (where
we do not require that $\mathcal{P}$ avoid $z$) and $\chi$ is the
{\em Euler measure}. Let $\mathcal{D}$ be a component of
$(-\Sigma')-\bigcup_i \alpha_i -\bigcup_i\gamma_i$.  Then
$\overline{n}_p(\mathcal{D})$ equals (i) $1$ if $p$ is in the
interior of $\mathcal{D}$, (ii) $0$ if $p$ does not intersect
$\mathcal{D}$, (iii) ${1\over 2}$ if $p$ is on an edge of
$\mathcal{D}$ (but not a corner), and (iv) ${1\over 4}$ if $p$ is on
a corner of $\mathcal{D}$. We then extend $\overline{n}_p$ linearly
to $\mathcal{P}$.

In the case at hand, the possible $x$'s and $y$'s are either in the
interior of $F$ or not in $F$, and therefore they either contribute
$1$ or $0$.  On the other hand, $\overline{n}_{z}(\mathcal{P})=-2$,
$\chi(\mathcal{P})=-2g(S)$, and
$\overline{n}_{z_0'}(\mathcal{P})=-1$ are constant, and it follows
that $\langle c_1(\mathbf{s}_{\mathbf{y}'}),[\hat F]\rangle=2-2g(S)$
is the minimal value and it is attained when both $x$ and $y$ are
not in $F$. (In fact, $(\{\beta_0\}\cup\beta,\{\alpha_0\}\cup\alpha,
z,w)$ is a ``sutured Heegaard diagram'' in the sense of \cite{Ni}.)

\begin{figure}[ht]
\s
\begin{overpic}[width=4in]{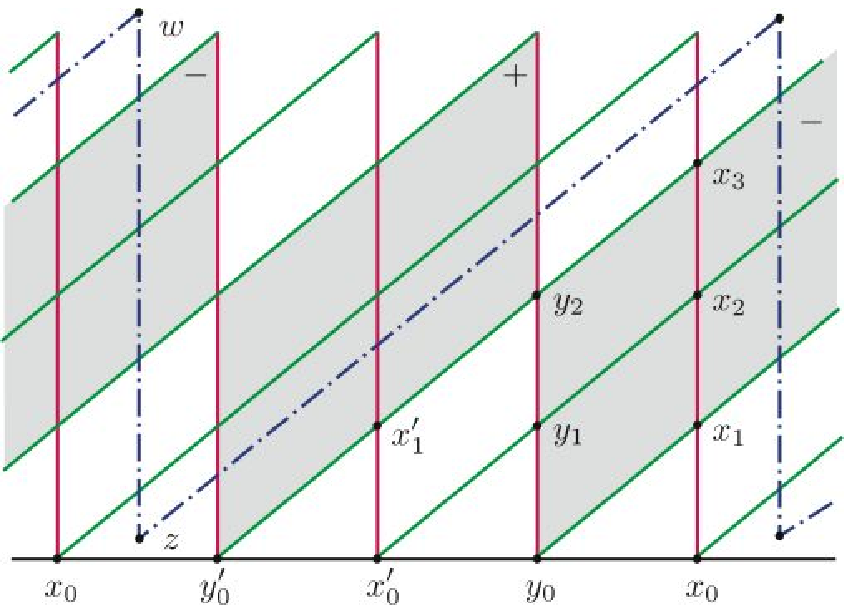}
\end{overpic}
\s
\caption{Cover of a neighborhood of $\bdry S$.} \label{spinc}
\end{figure}

The graded complex for calculating $\widehat{HFK}(-M,K,-2)$ is
generated by: $$(z_0, x_0,y_0), (z_0, x_0,y_2), (z_0, x_1',y_1),
(z_0, x_1,y_1), (z_0, x_2,y_0), (z_0, x_2,y_2), (z_0, x_3,y_1).$$
Our task is to identify $\mathbf{x}'=(z_0,x_0,y_0)$ as a generator
of $\widehat{HFK}(-M,K,-2)\simeq \Z$.  We will show that all the
generators besides $\mathbf{x}'$ correspond to \spinc structures
which are different from that of the contact structure $\xi$. An
easy computation shows that $H_2(M;\Z)\simeq \Z^2$ and is generated
by tori $T_\delta$ of the form $(\delta\times[0,1])/\sim$, where
$\delta$ is any nonseparating curve on $S$ and $(x,1)\sim (h(x),0)$
as before. Since $\xi$ is close to the foliation $S\times\{t\}$ on
$(S\times[0,1])/\sim$, it follows that $\langle c(\xi),[T_\delta]
\rangle=0$.  Now, let $\delta_1$ be a $(0,1)$-curve on $S$ and
$\delta_2$ be a $(1,0)$-curve.  Then $[T_{\delta_1}]$ is given by
the periodic domain $\mathcal{P}_{\delta_1}$, which consists of two
rectangles $y_0y_2y_4'y_2'$ and $y_0'y_2'y_4y_2$ with opposite
signs, shown in Figure~\ref{spinc}. Similarly, $[T_{\delta_2}]$ is
represented by $\mathcal{P}_{\delta_2}$, consisting of
$x_0x_2x_4'x_2'$ and $x_0'x_2'x_4x_2$ with opposite signs.

Now refer to Figure~\ref{spinc}, which is a cover of an annular
neighborhood of $\bdry S\subset S$. The dotted curve is (a lift of)
$c_1c_2$.  Points below the dotted curve are are not in $F$, so only
they have the proper filtration level to represent generators.

We will show $\mathbf{s}_{(z_0,x,y)}\not =
\mathbf{s}_{(z_0,x_0,y_0)}$ if $(x,y)\neq(x_0,y_0)$, by showing that
$\langle c(\mathbf{s}_{(z_0,x,y)}),[T_{\delta_i}] \rangle \not =0$
for $i=1$ and $2$ if $(x,y)\not=(x_0,y_0)$.

First consider the intersection points on the vertical lines
starting at $x_0$ and at $x_0'$. Suppose that $\langle
c(\mathbf{s}_{(z_0,x,y)}),[T_{\delta_2}] \rangle = 0$. The rectangle
$x_0'x_2'x_4x_2$ of the periodic domain $\mathcal{P}_{\delta_2}$
contributes $\frac{1}{2}$ if $x=x_3$ or $x=x_1'$. Since there is no
value of $y$ below the dotted curve with a contribution of
$-\frac{1}{2}$ from the rectangle $x_0x_2x_4'x_2'$ to cancel the
$1\over 2$, the possibilities $x=x_3,x_1'$ are eliminated. Since
$x_0x_2x_4'x_2'$ gives a contribution of $-{1\over 2}$ to $x_1$, and
$x_0'x_2'x_4x_2$ contributes $0$ to $y_0$, ${1\over 2}$ to $y_1$ and
$1$ to $y_1$, the only generator containing $x_1$ that is allowed is
$(z_0,x_1,y_1)$. Any generator containing $y_2$ is also disallowed
since $x_0'x_2'x_4x_2$ contributes $1$ to $y_2$, and there is no $x$
value that will offset it from the $x_0x_2x_4'x_2'$ rectangle. The
only generator allowed to contain $y_1$ is again $(z_0,x_1,y_1)$.
The same rectangle gives $x_1'$ a contribution of $-\frac{1}{2}$
that cannot be offset.

It therefore remains to consider the generator $(z_0,x_1,y_1)$, as
well as pairs with $x=x_0$ or $x_2$. Moreover, the only possible
$y$-coordinates are $y_0$ and $y_1$, and $(z_0,x_1,y_1)$ is the only
option allowed for $y=y_1$. Now use the periodic domain
$\mathcal{P}_{\delta_1}$. The rectangle $y_0y_2y_4'y_2'$ contributes
$-1$ to $(z_0,x_1,y_1)$, thus eliminating it as a possibility. The
only other option different from $(z_0,x_0,y_0)$ is $(z_0,x_2, y_0)$
(since $y_2$ was banned) which gets a nonzero contribution from
$y_0y_2y_4'y_2'$.

To show how this argument generalizes to higher genus surfaces, let
us examine the genus two case. The generators will have the form
$(z_0,x,y,u,v)$, and there will be 8 intersection points on each
vertical segment in a picture analogous to Figure~\ref{spinc}.
Denote the points on the boundary $u_0,v_0, u_0',v_0',x_0,y_0,
x_0',y_0'$ going from right to left. Start by considering the
rectangles $u_0u_2u_8'u_6'$ and $u_0'u_6'u_8u_2$. We eliminate $u_3,
\dots, u_7$ and all the $u'$ values besides $u'_0$, by noticing that
there is no allowable $v$ value to offset the $\frac{1}{2}$
contribution from $u_0'u_6'u_8u_2$. The contribution of $1$ from the
same rectangle eliminates all values of $v$ other than $v_0$ and
$v_1$ (though no $v_i'$ are yet disallowed). If $v=v_1$, only
generators of the form $(z_0,x,y, u_1,v_1)$ are allowed.

Now use the periodic domain represented by the rectangles
$v_0v_2v_8'v_6'$ and $v_0'v_6'v_8v_2$. The generators of the form
$(z_0,x,y, u_1,v_1)$ get a contribution of $-1$ from
$v_0v_2v_8'v_6'$ and there is no positive contribution from the
allowable $x,y$ coordinates that can be gained from
$v_0'v_6'v_8v_2$; therefore all the $(z_0,x,y,u_1,v_1)$ are
eliminated. Next, $u_2$ gets a contribution of $-1$ from
$v_0v_2v_8'v_6'$ that cannot be canceled since there is no $v$ value
that gets a contribution of $1$ needed from $v_0'v_6'v_8v_2$.
It follows that $u_0$ is the only allowable $u$-coordinate.
Generators $(z_0,x,y,u_0,v_i')$, $i\not=0$, are eliminated since
$v_i'$ gets a contribution of $\frac{1}{2}$ from $v_0'v_6'v_8v_2$
that cannot be canceled. Therefore we are left with $(z_0,x,y,
u_0,v_0)$.  The argument is now reduced to eliminating the possible
$x,y$ coordinates, and this follows just as in the genus one
argument given above.

This shows how the proof works for arbitrary genus. The inductive
step is done in the same way by eliminating all extra options in the
two new coordinates, thus reducing to the case of lower genus.

Since the contact invariant is the image of the generator of
$\widehat{HFK}(-M,K,-2g)$ in $\widehat{HF}(-M)$, it follows that
$c(\xi_{(S,h)})=EH(S,h)$. It is not hard to see how the above
argument generalizes to the  $h=\phi_{\bdry S}^n$, $n>0$ case.

\subsection{Isotopy} \label{isotopy}

In this subsection we prove the following:

\begin{lemma} \label{inv-isotopy}
If $h_t:S\stackrel\sim\rightarrow S$, $t\in[0,1]$, is a 1-parameter
family of diffeomorphisms which restrict to the identity on $\bdry
S$, then $EH(S,h_0,\{a_1,\dots,a_r\})= EH(S,h_1,\{a_1,\dots,a_r\})$.
\end{lemma}

\begin{proof}
Let $\alpha_i=(a_i,a_i)$ and $\beta_i^t=(b_i,h_t(b_i))$. In other
words, we fix the $\alpha_i$ and isotop the $\beta_i$.  Observe that
the $\beta_i^t$ remain constant on $S\times\{1\}$. According to
Theorem~7.3 of \cite{OS1}, we can reduce to the case where $h_t$ is
a Hamiltonian isotopy. Let $\Psi_t:\Sigma\stackrel\sim\rightarrow
\Sigma$ be the Hamiltonian isotopy which restricts to the identity
on $S\times\{1\}$ and restricts to $h_t$ on $S\times\{0\}$.  We use
the same notation for the induced isotopy on $Sym^r(\Sigma)$. Then
the chain map $\Phi:\widehat{CF}(\beta^0,\alpha)\rightarrow
\widehat{CF}(\beta^1,\alpha)$ is obtained by counting holomorphic
disks $u:[0,1]\times \R\rightarrow Sym^r(\Sigma)$ which satisfy
$\lim_{t\rightarrow +\infty} u(s+it)=\mathbf{x}$,
$\lim_{t\rightarrow -\infty} u(s+it)=\mathbf{x'}$, $u(0+it)\in
\Psi_t(\mathbb{T}_\beta)$, and $u(1+it)\in \mathbb{T}_\alpha$, and
avoid $\{z\}\times Sym^{r-1}(\Sigma)$. Here $\mathbf{x}\in
\widehat{CF}(\beta^0,\alpha)$ and $\mathbf{x'}\in
\widehat{CF}(\beta^1,\alpha)$. Now, if $\mathbf{x}$ is unique
$r$-tuple of points on $S\times\{1\}$ representing the generator of
$EH(S,h_0,\{a_1,\dots,a_r\}$, then the only holomorphic disk of the
above type are constant holomorphic disks, due to the placement of
the basepoint $z$.  This implies that $EH(S,h_0,\{a_1,\dots,a_r\})$
is mapped to $EH(S,h_1,\{a_1,\dots,a_r\})$ under the isomorphism
$\Phi: \widehat{HF}(\beta^0,\alpha)\stackrel\sim\rightarrow
\widehat{HF}(\beta^1,\alpha)$.
\end{proof}

\subsection{Change of basis} \label{change of basis}

In this subsection we prove the following proposition:

\begin{prop}\label{basischange}
$EH(S,h,\{a_1,\dots,a_r\})$ is independent of the choice of basis
$\{a_1,\dots,a_r\}$.
\end{prop}

Let $\{a_1,a_2,\dots,a_r\}$ be a basis for $S$. After possibly
reordering the $a_i$'s, suppose $a_1$ and $a_2$ are adjacent arcs on
$\bdry S$, i.e., there is an arc $\tau\subset \bdry S$ with
endpoints on $a_1$ and $a_2$ such that $\tau$ does not intersect any
$a_i$ in $\mbox{int}(\tau)$.  Define $a_1+a_2$ to be the isotopy
class of $a_1\cup\tau\cup a_2$, relative to the endpoints.  Then the
modification $\{a_1,a_2,\dots,a_r\}\mapsto
\{a_1+a_2,a_2,\dots,a_r\}$ is called an {\em arc slide}.

Proposition~\ref{basischange} is immediate from the following two
lemmas.

\begin{lemma}  \label{independence}
$EH(S,h)$ is invariant under an arc slide
$\{a_1,a_2,\dots,a_r\}\mapsto \{a_1+a_2,a_2,\dots,a_r\}$.
\end{lemma}

\begin{figure}[ht]
\begin{overpic}[width=2.6in]{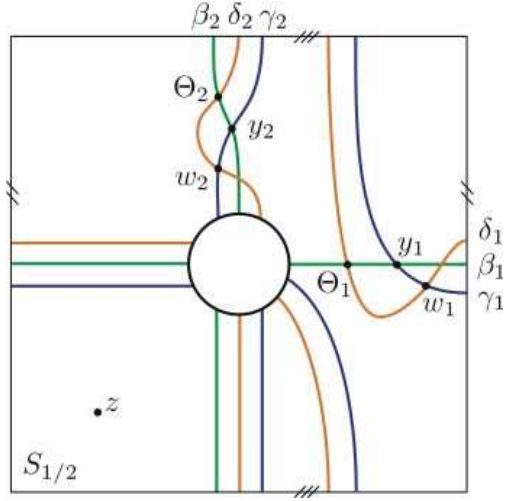}
\end{overpic}
\caption{The second handleslide.} \label{handleslide2}
\end{figure}

\begin{proof}
Without loss of generality, consider the case where $S$ is a
once-punctured torus.  We show that the chain map which corresponds
to an arc slide takes the representative of $EH(S,h, \{a_1,a_2\})$
determined by $\mathbf{x}=(x_1,x_2)$ to the representative of
$EH(S,h, \{a_1+a_2, a_2\})$ determined by the intersection point
$\mathbf{w}=(w_1,w_2)$. Observe that an arc slide corresponds to a
sequence of two handleslides for the corresponding Heegaard
splitting.

Let $(\Sigma, \beta, \alpha, z)$ be the pointed Heegaard diagram
corresponding to $a_i$, $b_i$ as described above, with $z$ a point
in $S_{1/2}$ lying outside the thin strips of isotopy between
$a_i$'s and $b_i$'s.  If we slide $\alpha_2$ over $\alpha_1$ along a
path parallel to $\bdry S$, then we obtain a new pair
$\gamma=\{\gamma_1,\gamma_2\}$, where $\gamma_1=(a_1+a_2,a_1+a_2)$
and $\gamma_2$ is a suitable pushoff of $(a_2, a_2)$ as in the proof
of the invariance of Heegaard Floer homology under handleslides in
\cite{OS1}. Figure~\ref{handleslide1} depicts the case where $a_1$
is to the right of $a_2$ with respect to $\tau$; the case where
$a_2$ is to the right of $a_1$ is treated similarly.

We claim that $(\Sigma,\gamma,\beta,\alpha,z)$ is a weakly
admissible Heegaard triple-diagram.  Recall that a triple-diagram is
{\em weakly admissible} if each nontrivial triply-periodic domain
which can be written as a sum of doubly-periodic domains has both
positive and negative coefficients.  First let us restrict to a
neighborhood $\mathcal{R}$ of the labeled regions of
$\Sigma-\cup_i\alpha_i-\cup_i\beta_i-\cup_i\gamma_i$ on the
right-hand side of Figure~\ref{handleslide1}. Due to the placement
of $z$, the only potential doubly-periodic region involving
$\beta,\alpha$ on $\mathcal{R}$ is $D_2+D_3-D_5-D_6$. (Here $D_i$ is
the domain labeled $i$.) Similarly, for $\gamma,\beta$ we have
$D_1+D_2-D_4-D_5$ and for $\alpha,\gamma$ we have $D_1+D_6-D_3-D_4$.
Taking linear combinations, we have
\begin{eqnarray*}
a(D_2+D_3-D_5-D_6)+b(D_1+D_2-D_4-D_5)+c(D_1+D_6-D_3-D_4) \\
=  (b+c)D_1+ (a+b)D_2+ (a-c)D_3 -(b+c)D_4 -(a+b)D_5 +(-a+c) D_6.
\end{eqnarray*}
Since the coefficients come in pairs, e.g., $a+b$ and $-(a+b)$, if
any of $a+b$, $b+c$, $a-c$ does not vanish, then the triply-periodic
domain has both positive and negative coefficients. Hence, if any of
$\alpha_1$, $\beta_1$ and $\gamma_1$ is used, then we are done.
Otherwise, we may assume that none of $\alpha_1$, $\beta_1$ and
$\gamma_1$ is used in the periodic domain.  This allows us to erase
all three, and apply the above considerations to $\alpha_2$,
$\beta_2$, and $\gamma_2$.  The verifications of weak admissibility
of all other triple-diagrams in this paper are identical, and are
omitted.

Let $\Theta=(\Theta_1,\Theta_2)$ be the top generator of
$\widehat{HF}(\# (S^1\times S^2))= \widehat{HF}(\alpha,\gamma)$.
Define the map
$$\psi:\widehat{HF}(\beta,\alpha)\otimes
\widehat{HF}(\alpha,\gamma)\rightarrow \widehat{HF}(\beta,\gamma),$$
where $\psi(\mathbf{y}\otimes \mathbf{y}')$ counts holomorphic
triangles, two of whose vertices are $\mathbf{y}\in
\widehat{CF}(\beta,\alpha)$ and $\mathbf{y}'\in
\widehat{CF}(\alpha,\gamma)$. Then the isomorphism
$g:\widehat{HF}(\beta,\alpha)\stackrel\sim\rightarrow
\widehat{HF}(\beta,\gamma)$ is given by
$g(\mathbf{y})=\psi(\mathbf{y}\otimes \Theta)$.

\begin{figure}[ht]
\begin{overpic}[width=2.6in]{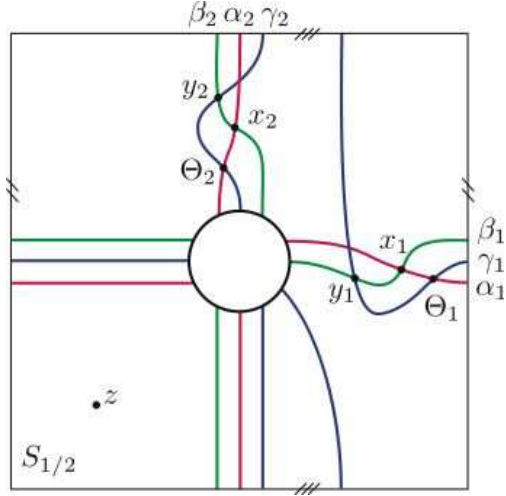}
\end{overpic}
\caption{The first handleslide.} \label{handleslide1}
\end{figure}

We claim that the representative $\mathbf{x} =(x_1,x_2)$ of
$EH(S,h,\{a_1,a_2\})$ gets mapped to $\mathbf{y}=(y_1,y_2)\in
\widehat{CF}(\beta,\gamma)$ given in Figure~\ref{handleslide1}.  By
the placement of $z$, we see that the unique holomorphic map $\hat
u$ which has $x_1$ and some $\Theta_i$ as corners (and avoids $z$)
must be a triangle with vertices $x_1,\Theta_1,y_1$.  Now that
$\alpha_1$, $\beta_1$, and $\gamma_1$ are used up, it easily follows
that the unique holomorphic map $\hat u$ which involves $x_2$ and
$\Theta_2$ (and avoids $z$) is a triangle with vertices
$x_2,\Theta_2,y_2$. This proves the claim.

Let us now consider the effect of the second handleslide, depicted
in Figure~\ref{handleslide2}.  Let $\delta=\{\delta_1,\delta_2\}$,
where $\delta_1$ and $\delta_2$ are suitable pushoffs of $(a_1+a_2,
h(a_1+a_2))$ and $(a_2, h(a_2))$, respectively. A similar argument
as above shows that, under the map
$$\widehat{HF}(\delta,\beta)\otimes
\widehat{HF}(\beta,\gamma)\rightarrow \widehat{HF}(\delta,\gamma),$$
$\Theta\otimes\mathbf{y}$ gets mapped to $\mathbf{w}$. This shows
that $\bf x$ and $\bf w$ determine the same element in Heegaard
Floer homology, and consequently
$EH(S,h,\{a_1,a_2\})=EH(S,h,\{a_1+a_2,a_2\})$.
\end{proof}

\begin{lemma}\label{handleslidesexist}
Let $\{a_1,\dots,a_{r}\}$ and $\{b_1,\dots,b_{r}\}$ be two bases for
$S$.  Then there is a sequence of arc slides that takes
$\{a_1,\dots,a_{r}\}$ to $\{b_1,\dots,b_{r}\}$.
\end{lemma}

We do not need to assume that $\bdry S$ is connected.

\begin{proof}
We argue that we can reduce the total number of intersections of
$\bigcup_i a_i$ and $\bigcup_i b_i$ by replacing
$\{a_1,\dots,a_{r}\}$ with $\{a_1',\dots,a_{r}'\}$, which is
obtained from $\{a_1,\dots,a_r\}$ by a sequence of arc slides.  By
inducting on the number of intersection points, this shows that we
can perform a sequence of arc slides until $\bigcup_i a_i$ and
$\bigcup_i b_i$ become disjoint.  We then show that two disjoint
bases can be brought one into another by a sequence of arc slides.

Let $P=S-\bigcup_i a_i$.  Then $P$ is a polygon whose boundary
$\bdry P$ consists of $4r$ arcs, $2r$ of which are $a_i$ or
$a_i^{-1}$ and $2r$ of which are arcs $\tau_1,\dots,\tau_r$ of
$\bdry S$.

Suppose $(\bigcup_i a_i)\cap (\bigcup_i b_i)\not=\emptyset$, where
we are assuming efficient intersections.  After possibly reordering
the arcs, there is a subarc $b_1^0\subset b_1$ which starts on
$\tau_1\subset \bdry S$ and ends on $a_1$, and whose interior
$int(b_1^0)$ does not intersect $\bigcup_i a_i$.  (In other words,
$b_1^0$ is a properly embedded arc of $P$.)  We may assume that
$a_1$ is not adjacent to $\tau_1$; otherwise, isotop the relevant
endpoint of $b_1$ along $\tau_1$. The subarc $b_1^0$ separates the
polygon $P$ into two regions $P_1$ and $P_2$, only one of which
contains a boundary arc that is labeled $a_1^{-1}$ (say $P_2$). We
can then slide $a_1$ over all the arcs of type $a_i$ or $a_i^{-1}$
in the other region $P_1$, and obtain the new curve $a_1'$ as in
Figure~\ref{handleslidesexistfig} so that the new basis
$\{a_1',a_2,\dots,a_r\}$ has fewer intersections with $\bigcup_i
b_i$.
\begin{figure}[ht]
\begin{overpic}[width=8cm]{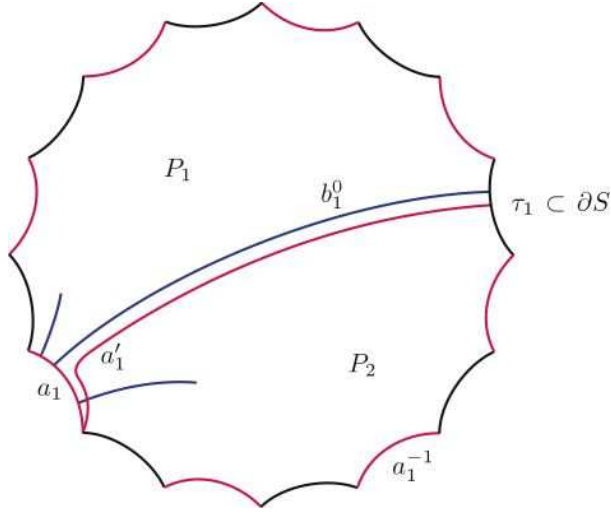}
\end{overpic}
\caption{Simplifying the intersections of $\bigcup a_i$ and $\bigcup
b_i$.} \label{handleslidesexistfig}
\end{figure}
(Note that trying to slide over $a_1^{-1}$ presents a problem, so we
must go the other way around.) There is one situation when the above
strategy needs a little more thought, namely when $\bdry P_2$ only
intersects $a_1$ and $a_1^{-1}$ (among all the $a_i$ and
$a_i^{-1}$).  In this case, $b_1$ exits the polygon $P$ along $a_1$
and reenters through $a_1^{-1}$.  Eventually we find a subarc of
$b_1$ which starts on some $\tau_2$ and ends on an adjacent
$a_1^{-1}$, a contradiction.  We now apply the same procedure to
$\{a_1',a_2,\dots,a_r\}$ and $\{b_1,\dots,b_r\}$ until they become
disjoint.

Now suppose that the two bases $\{a_1,\dots,a_r\}$ and
$\{b_1,\dots,b_r\}$ are disjoint.  We consider the polygon
$P=S-\bigcup_i a_i$. Some of the $b_i$ arcs may be parallel to $a_j$
or $a_j^{-1}$. An arc $b_1$ that is not parallel to any of the $a_i$
will cut $P$ into two components $P_1$ and $P_2$, each containing
more than one of $a_i, a_i^{-1}$, $i=1,\dots,r$. Recall that $b_1$
is nonseparating. One can easily verify that $b_1$ being
nonseparating is equivalent to the existence of some $a_i$ such that
$a_i\in P_1$ and $a_i^{-1}\in P_2$ (or vice versa). (If there is
some $a_i$, then take an arc $c$ in $P$ from $a_i\subset P_1$ to
$a_i^{-1}\subset P_2$. The closed curve in $S$ obtained by gluing up
$c$ is dual to $b_1$.) If each such $a_i$ is parallel to some $b_j$,
then $S-\bigcup_i b_i$ would be disconnected, so we could
additionally assume that there is some $a_i$ which is not parallel
to any $b_j$. Now we slide $a_i$ across all the arcs of type $a_j$,
$a_j^{-1}$ in $P_1$ until it becomes parallel to $b_1$.
\end{proof}

\subsection{Legendrian surgery} \label{Legendrian}

Let $\delta$ be a nonseparating curve and $\phi_\delta^{-1}$ be a
negative Dehn twist about $\delta$.   We now transfer $EH$ from
$M=M_{(S,h)}$ to $M'=M_{(S,\phi_\delta^{-1}\circ h)}$. Recall that
there is a natural map $$f:\widehat{HF}(-M)\rightarrow
\widehat{HF}(-M'),$$ which arises from tensoring with the top
generator $\Theta$ of $\widehat{HF}(\#(S^1\times S^2))$.

\begin{prop}
$f(EH(S,h))=EH(S,\phi_\delta^{-1}\circ h)$.
\end{prop}

\begin{proof}
By Proposition~\ref{basischange} we may take a basis
$\{a_1,\dots,a_r\}$ for $S$ so that $\delta$ is disjoint from
$h(b_2),\dots,h(b_r)$, intersects $h(b_1)$ exactly once, and is
parallel to $h(b_2)$. Then the result of performing $(+1)$-surgery
along $\delta$ (or, equivalently, a negative Dehn twist along
$\delta$) is given by Figure~\ref{surgery}.

\begin{figure}[ht]
\begin{overpic}[width=5in]{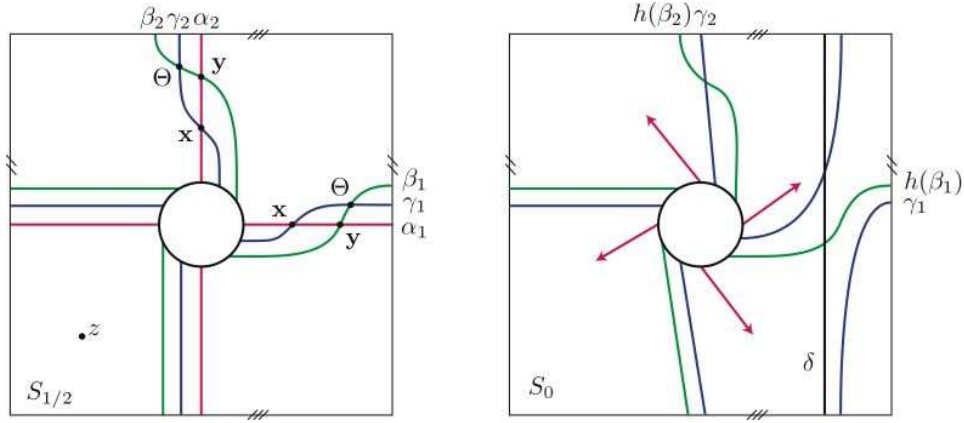}
\end{overpic}
\caption{Legendrian $(+1)$-surgery. The second figure shows curves
and their {\it relative} position correctly, but they are positioned
on $S_0$ as if $h^{-1}$ had been applied to each of them.}
\label{surgery}
\end{figure}

The $\alpha$-curves and $\beta$-curves are as before, and we define
the $\gamma$-curves as follows: Let
$\gamma_1=(b_1,\phi_{\delta}^{-1}\circ h(b_1))$ and $\gamma_i=(b_2,
h(b_i))$ for $i>1$. Let $\Theta\in \widehat{HF}(\gamma,\beta)$ be
the top generator of $\# (S^1\times S^2)$, given in
Figure~\ref{surgery}. Define the map
$$\phi:\widehat{HF}(\gamma,\beta)\otimes
\widehat{HF}(\beta,\alpha)\rightarrow \widehat{HF}(\gamma,\alpha),$$
where $\phi(\mathbf{y}\otimes\mathbf{y'})$ counts holomorphic
triangles, two of whose vertices are $\mathbf{y}$ and $\mathbf{y}'$.
Then the map $f:\widehat{HF}(\beta,\alpha)\rightarrow
\widehat{HF}(\gamma,\alpha)$ is given by
$f(\mathbf{y})=\phi(\Theta\otimes \mathbf{y})$.  By the convenient
placement of $z$, it follows that we only have small triangles in
the Heegaard diagram. Hence if $[\mathbf{x}]= EH(S,h,\{a_1,a_2\})$,
then $\phi([\Theta\otimes \mathbf{x}])= EH(S,\phi_\delta^{-1}\circ
h,\{a_1,a_2\})$.
\end{proof}

\subsection{Multiple boundary components}\label{multiple}
Consider $(S,h)$ where $S$ has disconnected boundary.  For
simplicity, assume $S$ has two boundary components. Pick a basis
$\{a_1,\dots,a_r\}$ for $S$.  Next consider $(S',h\# id)$, where
$S'$ is obtained from $S$ by attaching a 1-handle between the two
boundary components and we are extending $h$ by the identity. If
$a_0$ is the cocore of the 1-handle, then $\{a_0,\dots,a_r\}$ is a
basis for $S'$.  Our argument is similar to that of Lemma~4.4 of
\cite{OS3}. The natural map
$$F_U:\widehat{HF}((-M_{(S,h)})\#(S^1\times S^2))\rightarrow \widehat{HF}(-M_{(S,h)}),$$
which corresponds to the cobordism $U$ attaching a 3-handle as in
Section~4.3 of \cite{OS4}, sends
$$EH(S',h\# id,\{a_0,a_1,\dots,a_r\})\mapsto EH(S,h,\{a_1,\dots,a_r\}).$$
Since $S'$ has only one boundary component, we already know that
$$c(S',h\# id)=EH(S',h\# id,\{a_0,a_1,\dots,a_r\}).$$ Moreover, if
$\delta$ is a closed curve on $S'$ which is ``dual'' to $a_0$, then
there is a natural map
$$F_W: \widehat{HF}(-M_{(S,h)})\rightarrow \widehat{HF}((-M_{(S,h)})\#(S^1\times
S^2))$$ which maps $c(S,h)$ to $c(S',h\# id)$.  Here $(S,h)$ and
$(S', \phi_\delta \circ(h\# id))$ represent the same 3-manifold, and
$W$ is the cobordism corresponding to the Legendrian $(+1)$-surgery.
Finally, $U\circ W\simeq [0,1]\times M_{(S,h)}$, so
$$c(S,h)=F_U\circ F_W(c(S,h))=F_U(c(S',h\#
id))=EH(S,h,\{a_1,\dots,a_r\}).$$

\section{Right-veering and holomorphic disks} \label{EH}

In this section we prove Theorem~\ref{torus}.

\begin{proof}[Proof of Theorem~\ref{torus}.]
Let $S$ be a once-punctured torus.

Suppose first that $h$ has pseudo-Anosov monodromy.  If the
fractional Dehn twist coefficient $c\geq 1$, then the contact
structure is already symplectically fillable and universally tight.
It also follows that $c(\xi_{(S,h)})\not=0$.  If $c=\frac{1}{2}$,
then $c(\xi_{(S,h)})\not=0$ follows from Theorem~\ref{non-zero}
below. If $c\leq 0$, then $\xi$ is overtwisted since $S$ is not
right-veering. (See \cite{HKM2}.)

If $h$ is periodic, then $\xi$ is right-veering if and only if $h$
is a product of positive Dehn twists by \cite{HKM3}.

If $h$ is reducible, then $c(\xi_{(S,h)})\not=0$ follows from
Theorem~\ref{nonzero2} below.
\end{proof}

\begin{thm}\label{non-zero}
Let $(S,h)$ be an open book decomposition for $M$, where $S$ is a
once-punctured torus and $h$ is pseudo-Anosov with fractional Dehn
twist coefficient $c = \frac{1}{ 2}$.  Then $c(\xi_{(S,h)}) =EH(S,h)
\neq 0$, and hence the contact structure $\xi_{(S,h)}$ is tight.
\end{thm}

\begin{proof}
We show that $EH(S,h) \neq 0$ by choosing a basis for $S$ for which
there are no holomorphic disks in the corresponding Heegaard diagram
that map to the generator $\mathbf{x}=(x_0,y_0)$ defining $EH(S,h)$.

The following lemma furnishes us with a convenient basis:

\begin{lemma}
Let $A \in SL(2,\mathbb Z)$ be a matrix with $tr(A)< -2$. Then $A$
is conjugate in $SL(2, \mathbb Z)$ to a matrix $\begin{pmatrix} a & b\\
c & d\end{pmatrix}$, where $(a,c)$ and $(b,d)$ are in the third
quadrant.
\end{lemma}

\begin{proof}
Let $\Lambda^s$ and $\Lambda^u$ be the stable and unstable
laminations for $A$. The slopes of $\Lambda^s$ and $\Lambda^u$ will
be written $\mbox{slope}(\Lambda^s)$ and $\mbox{slope}(\Lambda^u)$.
(Recall that these slopes are irrational.) Let us consider the Farey
tessellation on the hyperbolic unit disk $D^2$.  Pick a vertex $s_1$
on the clockwise edge along $\bdry D^2$ from
$\mbox{slope}(\Lambda^s)$ to $\mbox{slope}(\Lambda^u)$, and pick a
vertex $s_2$ on the counterclockwise edge from
$\mbox{slope}(\Lambda^s)$ to $\mbox{slope}(\Lambda^u)$, so that
there is an edge of the Farey tessellation between $s_1$ and $s_2$.
(The existence of such a pair $s_1,s_2$ is an exercise.) Then
$A(s_1)$ (resp.\ $A(s_2)$) is closer to $\mbox{slope}(\Lambda^s)$
than $s_1$ (resp.\ $s_2$) is. An oriented basis corresponding to
$(s_1,s_2)$ will have the desired property.
\end{proof}

With the choice of basis as above, we can represent $M=M_{(S,h)}$ by
the Heegaard diagram below. We have drawn a picture of the diagram
corresponding to $A=\begin{pmatrix} -1 & -1 \\ -1 & -2
\end{pmatrix}$, but the same argument works for any such $A$ as
described in the previous lemma.  We prove that there is no
holomorphic disk from any $\mathbf{y}$ to $\mathbf{x}=(x_0,y_0)$.
Suppose on the contrary that there is such a holomorphic disk $u$.
Let $\hat u:\widehat D\rightarrow \Sigma$ be the corresponding
holomorphic map to $\Sigma$.  Assuming $\bdry \widehat D$ is
connected, it is given by a subarc of $a_1$ from some $x_i\in
a_1\cap h(a_2)$ to $x_0$, followed by a subarc of $h(a_1)$ from
$x_0$ to some $y_j\in a_2\cap h(a_1)$, followed by a subarc of $a_2$
from $y_j$ to $y_0$ (you either turn left or turn right at $y_j$),
and then by a subarc of $h(a_2)$ from $y_0$ to $x_i$. If we lift
$\bdry \widehat{D}$ to the universal cover of the capped off surface
$T^2=S\cup D^2$, then in all cases we see that $\bdry \widehat{D}$
is not contractible.  This implies that $\bdry \widehat{D}$ cannot
bound a surface in $S$.  We argue similarly when $\bdry \widehat D$
has two components.  It follows that the class $EH(S,h)$ of
$\mathbf{x}=(x_0,y_0)$ is nonzero.
\end{proof}

\begin{figure}[ht]
\begin{overpic}[width=12cm]{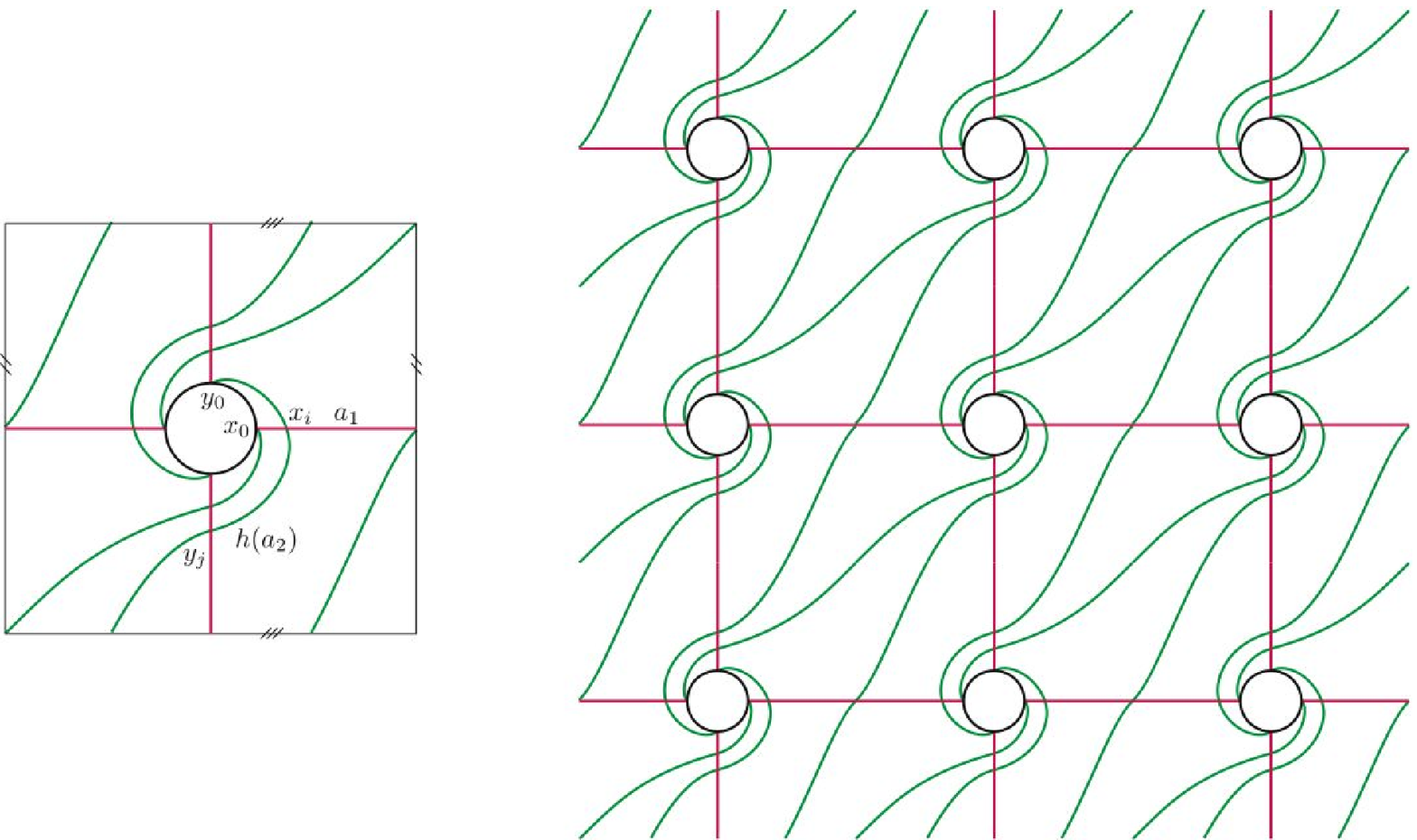}
\end{overpic}
\caption{} \label{puncturedtorus}
\end{figure}

\begin{thm}\label{nonzero2}
$EH(S,h)\not=0$ if $h$ is reducible and right-veering.
\end{thm}

\begin{proof}
Suppose $h$ is reducible.  Let $g$ be an element of $Aut(S,\bdry S)$
which is the minimally right-veering representative for the matrix
$A=-id$. (In terms of positive Dehn twists, $g=(A_1A_2A_1)^2$, where
$A_1=\begin{pmatrix} 1 & 1\\ 0 & 1\end{pmatrix}$ and
$A_2=\begin{pmatrix} 1 & 0 \\ -1 & 1 \end{pmatrix}$.)   After
changing bases if necessary, $h=g^n \phi_\gamma^m$, where $n$ is a
positive integer, $m$ is an integer, and $\phi_\gamma$ is a positive
Dehn twist about a $(0,1)$-curve $\gamma$. If $m$ is nonnegative,
then $h$ is a product of positive Dehn twists, and $EH(S,h)\not=0$.

Suppose $m<0$.  It suffices to prove the theorem for $n=1$, since
the contact structures corresponding to larger $n$ are obtained from
the $n=1$ case by Legendrian surgery. Take a basis corresponding to
slopes $0,\infty$ and matrix $A=\begin{pmatrix} -1 & 0 \\ -m &
-1\end{pmatrix}$.  Then $EH(S,h)$ is nonzero by the same method as
in Theorem~\ref{non-zero}.
\end{proof}

\s\n {\em Acknowledgements.}  The authors are grateful to John
Etnyre for discussions which led to the alternate description of the
contact class.

\end{document}